\documentclass[screen]{aomart}

\makeatletter
\fancypagestyle{firstpage}{%
  \fancyhf{}%
  \if@aom@manuscript@mode
    \lhead{\begin{picture}(0,0)%
        \put(-21,-25){\usebox{\@aom@linecount}}%
      \end{picture}}
  \fi
  \chead{}%
   \cfoot{\footnotesize\thepage}}%
\makeatother

\usepackage{cancel}
\usepackage{amsthm}
\usepackage{mathtools}
\chardef\bslash=`\\ 

\newtheorem[{}\it]{thm}{Theorem}
\newtheorem{cor}[thm]{Corollary}
\newtheorem{lem}[thm]{Lemma}
\newtheorem{prop}[thm]{Proposition}

\newtheorem{afa}{Fact}

\theoremstyle{definition}
\newtheorem{defn}{Definition}
\newtheorem{rem}{Remark}
\newtheorem{exa}{Example}

\newtheorem*[{}\it]{notation}{Notation}

\addresses

\title[Real Linear $2^\mathrm{nd}$ order Ordinary diff. eqs. and hyperbolic geometry]{Almost all real linear second order ordinary differential equations are solved by geodesic curves in two dimensional Riemannian hyperbolic geometry
}
\author{{\L}ukasz Rudnicki}
\address{Faculty of Mathematics, Physics, and Informatics \& International Centre for Theory of Quantum Technologies, 
 University of Gdansk,
80-308 Gdansk, Poland}
\email{lukasz.rudnicki@ug.edu.pl} 
\orcid{0000-0001-8563-6101}


\copyrightnote{}

\begin{document}
\begin{abstract}
I show that a real linear second order ordinary differential equation $u''\left(x\right)+h\left(x\right)u\left(x\right)=0$, with differentiable $h(x)$, locally admits two linearly independent solutions which exist on an open interval around any $x_0\in\mathbb{R}$:
\[
        u_\mathtt{top}(x)=\exp\left[\int_{x_0}^{x}\!\!d\xi\,\Phi\left(\xi\right)\frac{\Phi'\left(\xi\right)-\sqrt{\left[h\left(\xi\right)-\Phi^{2}\left(\xi\right)\right]^{2}+\left[\Phi'\left(\xi\right)\right]^{2}}}{h\left(\xi\right)-\Phi^{2}\left(\xi\right)}\right],
\]
\[
        u_\mathtt{bot}(x)=\exp\left[\int_{x_0}^{x}\!\!d\xi\,\Phi\left(\xi\right)\frac{\Phi'\left(\xi\right)+\sqrt{\left[h\left(\xi\right)-\Phi^{2}\left(\xi\right)\right]^{2}+\left[\Phi'\left(\xi\right)\right]^{2}}}{h\left(\xi\right)-\Phi^{2}\left(\xi\right)}\right],
\]
where $\Phi(x)$ is any geodesic curve in a two dimensional \textit{hyperbolic} geometry of a Riemannian manifold $\mathbb{M}_h$, which is non-vertical at $x_0$. I define  $\mathbb{M}_h$ to be  an upper half plane $\{\left(x,\varPhi\right)\in\mathbb{R}^2\,|\,\varPhi>0\}$, with points in which $\varPhi^2=h(x)$ being removed, equipped with metric $g_h=\left[\left(h(x)-\varPhi^2\right)^2dx^2+d\varPhi^2\right]/\varPhi^2$. A non-trivial character of the presented result stems from the fact that $g_h$ is solely defined in terms of the function $h(x)$.

I also show that a local diffeomorphism between $\mathbb{M}_h$ and Poincar{\'e} upper half plane $\mathbb{H}$ is induced by any pair of linearly independent solutions of $u''\!\left(x\right)+h\left(x\right)u\left(x\right)=0$. If this pair is selected to be $u_\mathtt{top}(x)$ and $u_\mathtt{bot}(x)$, the associate geodesic curve $\Phi(x)$ is mapped to a vertical geodesic curve on $\mathbb{H}$.

With these results, supported by complementary remarks and examples, I establish a fundamental link between linear second order ordinary differential equations and two dimensional hyperbolic geometry, adding to textbook knowledge about these profound subareas of mathematics.
\end{abstract}

\maketitle
\tableofcontents

\section{Introduction}
We start with a general linear second order ordinary differential equation 
\begin{equation}\label{Eq.1}
\left(a\left(t\right)u'\left(t\right)\right)'+b\left(t\right)u\left(t\right)=c\left(t\right).
\end{equation}
We assume the independent variable $t$ is real and the functions $a\left(t\right)$, $b\left(t\right)$ and $c\left(t\right)$ are real-valued, therefore, we develop a real theory. However, even if two special solutions of (\ref{Eq.1}) are taken as real, a general solution can be complex due to linearity. From now on we skip the adjective \textit{real}. The above three functions are defined on $\mathbb{R}$, or on an interval, though, for the sake of clarity we assume the former case (if need be, one can always restrict the domain).
\begin{defn}[Almost all linear-$2^{nd}$-order-ODEs]\label{Almost}
   The subject of our study are Eqs. (\ref{Eq.1}) with $a\left(t\right)$, $b\left(t\right)$ and $c\left(t\right)$ continuous (of class $C^0$), $a\left(t\right)\neq 0$ and $a\left(t\right)b\left(t\right)$ differentiable (of class $C^1$).
\end{defn}
Since continuity of all three functions and positivity or negativity of $a\left(t\right)$ are minimal assumptions to assure existence and uniqueness of initial value problem associated with (\ref{Eq.1}), the adjective \textit{almost} refers to an additional differentiability requirement imposed on the product $a\left(t\right)b\left(t\right)$.

 When $c\left(t\right)\equiv 0$, Eq. (\ref{Eq.1}) is called homogeneous.
If two linearly independent solutions of the homogeneous equation are known, the general solution of (\ref{Eq.1}) can be expressed by quadratures. Therefore,  without loss of generality we consider the homogeneous equation only. 

The homogeneous variant of Eq. (\ref{Eq.1}), subject to continuity assumptions collected in \fullref{Definition}{Almost}, can be reduced to the form \cite{Hartman}
\begin{equation}\label{Eq.main}
u''\left(x\right)+h\left(x\right)u\left(x\right)=0,
\end{equation}
where a differentiable function $h(x)$ and a solution $u(x)$ are defined as
\begin{equation}
    h(x) = a\left(t(x)\right)b\left(t(x)\right),\qquad u(x) = u\left(t(x)\right),
\end{equation}
and $t(x)$ is the inverse of $x(t)=\int_{t_0}^t\!\! d\tau\, a^{-1}(\tau)$, for some $t_0$. This transformation also holds without differentiability assumption, i.e. in general for all well-behaved equations, however, then $h(x)$ is only guaranteed to be of class $C^0$. For further convenience, let us now make a notation remark.
\begin{notation}[Solutions of Eq. (\ref{Eq.main})] With $u_1(x)$ and $u_2(x)$ we denote a pair of arbitrary (unspecified) linearly independent solutions of (\ref{Eq.main}). By $W$ we denote the Wronskian of these solutions. The most important result of this paper --- a pair of solutions which is locally represented in terms of a geodesic curve in two dimensional hyperbolic geometry --- will be denoted as  $u_\mathtt{top}(x)$ and $u_\mathtt{bot}(x)$. The subscripts "$\mathtt{top}/\mathtt{bot}$" refer to $\mathrm{"top"}/\mathrm{"bottom"}$. A reason for that choice is given below \fullref{Theorem}{thm-main2}. The Wronskian of these special solutions will be denoted as $W_\mathtt{t-b}$. Any specific pair of solutions and their Wronskian, relevant for a particular choice of $h(x)$, will be distinguished by an additional subscript following the number of the solution ($1$ or $2$), or the subscript "$\mathtt{top}/\mathtt{bot}$".
\end{notation}

The properties of Eq. (\ref{Eq.main}) and its solutions seem extremely deeply explored. Textbooks \cite{Hartman,Amann} rather concordantly cover: existence and uniqueness, the case of one special solution being known, reduction to Ricatti equation, Liouville substitution, Pr{\"u}fer transformations, Sturm comparison and separation theorems, Sturm Liouville theory, oscillation criteria, distributions of zeros, asymptotic expansions, monotones and more specialized results. 

In this contribution I show in \fullref{Sec.}{Sec2} that this equation, for any differentiable $h\left(x\right)$ can be solved by geodesic curves in \textit{the same} hyperbolic geometry in two dimensions, just expressed in different coordinates. In particular, these curves solve a geodesic equation on a Riemannian manifold $\mathbb{M}_h$, also introduced now, which is shown to be locally diffeomorphic to the paramount Poincar{\'e} upper half plane model of hyperbolic two dimensional geometry. The diffeomorphism can be built from any pair of linearly independent solutions of (\ref{Eq.main}), and in fact works for any continuous $h\left(x\right)$. The main result, presented as \fullref{Theorem}{thm-main2}, could be termed \textit{reduction of (\ref{Eq.main}) to the geodesic equation}, in a similar way one often refers to interrelation between (\ref{Eq.main}) and the Ricatti equation. However, despite an already quite non-trivial character of this assertion, the gist is in the surprising discovery that all $h$-dependent manifolds $\mathbb{M}_h$, on which the geodesic curves are to be found, always model the same, and at the same time very distinguished  hyperbolic geometry, and this is true \textit{without}  conformal equivalence assured by Uniformization Theorem.
\newpage
Already at the beginning it is good to stress that, solely due to the very nature of the presented results, the discussion is going to be conducted locally and in a coordinate-dependent fashion (e.g. referring to explicit Christoffel symbols). While differential geometry usually looks at coordinate-independent features of a manifold, we follow that path when showing that hyperbolic two dimensional geometry (a coordinate-independent feature) at once covers all differential equations specified by \fullref{Definition}{Almost}. This is true because the latter differ between each other exactly by what is a coordinate transformation between different representations of the same hyperbolic geometry.  

Following a supplementary discussion in \fullref{Sec.}{Sec3}, proofs of all results and details of presented examples are collected in \fullref{Sec.}{Sec4} and \fullref{Sec.}{Sec5} respectively.  

As the last step, before moving to the main part, I recall a definition of the Poincar{\'e} upper half plane.
\begin{defn}[Poincar{\'e} upper half plane]\label{DPoinc}
A set $\mathcal{H} = \{\left(X,Y\right)\in\mathbb{R}^2 \,|\, Y>0\}$ is called the upper half plane. It is a smooth manifold. Let 
\begin{equation}
    g_\mathbb{H}=\frac{dX^2+dY^2}{Y^2},
\end{equation}
be a Riemannian metric on $\mathcal{H}$. Poincar{\'e} upper half plane $\mathbb{H}=\left(\mathcal{H}, g_\mathbb{H}\right)$ is the manifold $\mathcal{H}$ equipped with the metric $ g_\mathbb{H}$. 
\end{defn}
The pair $\mathbb{H}=\left(\mathcal{H}, g_\mathbb{H}\right)$ is a complete simply connected Riemannian manifold with constant sectional curvature equal to $-1$.
\section{Main results}
\label{Sec2}
Let us now introduce a different geometry on the upper half plane $\mathcal{H}$.
\begin{defn}[Linear-$2^\mathrm{nd}$-order-ODE upper half plane] \label{DefM}
Let us express upper half plane in different coordinates as $\mathcal{H} = \{\left(x,\varPhi\right)\in\mathbb{R}^2\,|\,\varPhi>0\}$ and let
\begin{equation}\label{metric}
    g_h=\frac{\left(h(x)-\varPhi^2\right)^2dx^2+d\varPhi^2}{\varPhi^2},
\end{equation}
be a metric on $\mathcal{H}$. As $\det g_h=0$ if $\varPhi^2=h(x)$, this metric is degenerate. Therefore, let $\mathcal{H}_0 =\left\{\left(x,\varPhi\right)\in\mathcal{H}\,|\,\varPhi^2\neq h(x)\right\}$. A pair $\mathbb{M}_h = \left(\mathcal{H}_0, g_h\right)$ is a Riemannian manifold called the \textit{Linear-$2^\mathrm{nd}$-order-ODE upper half plane}.
\end{defn}
\begin{rem}
    Observe that $\mathbb{M}_h$ does not need to be connected. For certain choices of $h(x)$, for example $h(x)=\sin{x}$, this manifold is composed of infinitely many disjoint pieces.
\end{rem}
I call $\mathbb{M}_h$ "Linear-$2^\mathrm{nd}$-order-ODE upper half plane" despite the fact that a sole inclusion of an arbitrary function $h(x)$ does not yet bring a relation between $g_h$ and (\ref{Eq.main}). However, the presented results will justify that choice.
\subsection{Curvature of $g_h$}\label{subsecCur}
We start with:
\begin{lem}\label{lem-hyp}
Let $h(x)$ be of class $C^2$, so that the Riemann tensor of $g_h$ can be defined. The metric $g_h$ in the whole domain $\mathcal{H}_0$ has constant sectional curvature equal to $-1$. Therefore, it everywhere describes hyperbolic geometry.
\end{lem}
\begin{rem}\label{RemC0}
    In \fullref{Lemma}{lem-hyp} the function $h(x)$ is customarily assumed to be of class $C^2$, since the Riemann tensor involves second derivatives of the metric. However, in the case of $g_h$ all terms involving first and second derivatives of $h(x)$ disappear due to antisymmetry of the Riemann tensor. Therefore, treating antisymmetrization like a regularization procedure removing potential singularities of $h'(x)$ and $h''(x)$, we can extend the validity of \fullref{Lemma}{lem-hyp} to functions $h(x)$ which are absolutely continuous. Moreover, since the Riemann tensor of $g_h$ is regular when $\varPhi^2=h(x)$, the degeneracy of $g_h$ happening in this case stems from the choice of the coordinates. This is very similar to the degeneracy of polar coordinates in $\mathbb{R}^2$, which occurs in the origin. 
\end{rem}

As a consequence of \fullref{Lemma}{lem-hyp}, we know there shall exist a local diffeomorphism from $\mathbb{M}_h$ to $\mathbb{H}$. I will later prove (see \fullref{Lemma}{thm-main}) that such a diffeomorphism can be constructed from the solutions of  (\ref{Eq.main}). In relation to \fullref{Remark}{RemC0}, indeed it will turn out that the Jacobian of this map is continuous if $h\in C^0\left(\mathbb{R}\right)$, on the other hand, it vanishes for $\varPhi^2=h(x)$.

\subsection{Geodesic curves on $\mathbb{M}_h$}\label{subsecGeo} With the dot we  denote derivatives with respect to a parameter "$s$" while, as before, the prime is reserved for the derivative with respect to the variable "$x$". In local coordinates $\varsigma^k$, the geodesic equation \cite{Besse}  $\dot{\varsigma}^j\nabla_j\dot{\varsigma}^i=0$ reads $\ddot{\varsigma}^i+\Gamma_{\:\:\:jk}^{i}\left(\varsigma\right)\dot{\varsigma}^j\dot{\varsigma}^k=0$, where $\nabla$ denotes the unique Levi-Civita connection and Einstein summation convention is used. We proceed with a definition specific to $\mathbb{M}_h$, substituting Christoffel symbols $\Gamma_{\:\:\:jk}^{i}$ explicitly given in the proof of \fullref{Lemma}{lem-hyp}. 
\begin{defn}
    A curve $\mathbb{R}\ni s \mapsto \left(x(s),\varPhi(s)\right)\in \mathcal{H}_0$ parametrized by an affine parameter $s$ is a geodesic curve on $\mathbb{M}_h$ if it satisfies differential equations: \begin{subequations}\label{geode}
\begin{equation}\label{geode1}
\ddot{x}\left(s\right)=\frac{h'\left(x\left(s\right)\right)}{\varPhi^{2}\left(s\right)-h\left(x\left(s\right)\right)}\left[\dot{x}\left(s\right)\right]^{2}-2\frac{\varPhi^{2}\left(s\right)+h\left(x\left(s\right)\right)}{\varPhi\left(s\right)\left[\varPhi^{2}\left(s\right)-h\left(x\left(s\right)\right)\right]}\dot{x}\left(s\right)\dot{\varPhi}\left(s\right),
\end{equation}
\begin{equation}\label{geode2}
\ddot{\varPhi}\left(s\right)=\frac{\varPhi^{4}\left(s\right)-h^{2}\left(x\left(s\right)\right)}{\varPhi\left(s\right)}\left[\dot{x}\left(s\right)\right]^{2}+\frac{1}{\varPhi\left(s\right)}\left[\dot{\varPhi}\left(s\right)\right]^{2},
\end{equation}     \end{subequations}
with some generic initial conditions: $x(0)=x_0$, $\dot{x}(0)=\dot{x}_0$, $\varPhi(0)=\varPhi_0$, $\dot{\varPhi}(0)=\dot{\varPhi}_0$. We assume that $h(x)$ is of class $C^1$.
\end{defn}
From (\ref{geode1}) we can see that $\mathbb{M}_h$, similarly to $\mathbb{H}$, admits vertical geodesic curves $x(s)=x_0=\mathrm{const}$ and $\dot{x}(s)=\dot{x}_0=0$, which are of the form $\varPhi\left(s\right)=\varPhi_0 e^{\lambda s}$ where $\lambda=\dot{\varPhi}_0/\varPhi_0$. From this solution we conclude that if $\mathcal{H}\neq \mathcal{H}_0$, then $\mathbb{M}_h$ is not complete.

If within the initial conditions we chose $\dot{x}_0\neq 0$, a geodesic curve passing through $\left(x_0,\varPhi_0\right)$ is \textit{not vertical}. Standard existence theorems for ordinary differential equations assure \cite{Hartman} that for fixed initial conditions a unique geodesic curve of that kind exists in a neighborhood of $s=0$. In this neighborhood, let $s_-<0<s_+$ be such that 
\begin{equation}\label{inverdom}
    \displaystyle\mathop{\forall}_{s\in\left]s_{-},s_{+}\right[} \;\; \dot{x}\left(s\right)\neq 0.
\end{equation}
Then $x(s)$ is injective and consequently  invertible on $\left]s_{-},s_{+}\right[$. We define:
\begin{defn}[Geodesic curve on $\mathbb{M}_h$ in explicit form] If  a geodesic curve on $\mathbb{M}_h$, $\mathbb{R}\ni s \mapsto \left(x(s),\varPhi(s)\right)\in \mathcal{H}_0$, is such that the condition (\ref{inverdom}) holds, then we re-express the geodesic curve through a substitution $\varPhi(s)= \Phi\left(x(s)\right)$ and  call the function $\Phi\left(x\right)$ a geodesic curve on $\mathbb{M}_h$ represented in explicit form. It exists on $\left]x_-,x_+\right[$, where $x_-<x_0<x_+$. If $\dot{x}_0 > 0$, then  $x_-=x(s_-)$ and $x_+=x(s_+)$, while if $\dot{x}_0 < 0$, then  $x_-=x(s_+)$ and $x_+=x(s_-)$.
\end{defn}
\begin{rem}\label{hierar}
Since with the choice of the initial conditions we can assure that $\dot{x}_0\neq 0$, the above construction always works at least locally. It is also important to emphasize that the notion of the geodesic curve in explicit form is secondary with respect to solutions of (\ref{geode}). We first and foremost redefine $\varPhi(s)\coloneqq \Phi\left(x(s)\right)$ and then utilize invertibility of $x(s)$ to write $\Phi\left(x\right) = \varPhi\left(s(x)\right)$.
\end{rem}
\begin{notation}[Coordinates and geodesic curves in explicit form] With italic font I denote the coordinate $\varPhi$ and its $s$-dependent variant $\varPhi(s)$ (which is the coordinate of the geodesic curve), while functions of the variable $x$ representing geodesic curves in explicit form $\Phi(x)$ are denoted with upright $\Phi$.
\end{notation}
Before moving to the main result we need to establish the evolution equation for geodesic curves in explicit form. 
\begin{prop}\label{Prop1}
A geodesic curve on $\mathbb{M}_h$ in explicit form $\Phi(x)$ obeys
\begin{equation}\label{ELexp}
 \Phi''\left(x\right) 
 =\frac{3\Phi^{2}\left(x\right)+h\left(x\right)}{\Phi^{2}\left(x\right)-h\left(x\right)}\frac{\left[\Phi'\left(x\right)\right]^{2}}{\Phi\left(x\right)}-\frac{h'\left(x\right)\Phi'\left(x\right)}{\Phi^{2}\left(x\right)-h\left(x\right)}+\frac{\Phi^{4}\left(x\right)-h^{2}\left(x\right)}{\Phi\left(x\right)}.
\end{equation}
\end{prop}
We note in passing  that Eq. (\ref{ELexp}) can equivalently be written as
\begin{equation}
 \Phi''\left(x\right) =\Lambda\left(x,\Phi(x)\right)\left[\Phi'(x)\right]^{2}+\Gamma_{\:\:\:xx}^{x}\left(x,\Phi(x)\right)\Phi'(x)-\Gamma_{\:\:\:xx}^{\varPhi}\left(x,\Phi(x)\right),
\end{equation}
where $\Lambda\left(x,\Phi(x)\right)=2\Gamma_{\:\:\:\varPhi x}^{x}\left(x,\Phi(x)\right)-\Gamma_{\:\:\:\varPhi\varPhi}^{\varPhi}\left(x,\Phi(x)\right)$.

It is a peculiar feature of the geodesic equation, stemming from the fact that an acceleration depends only on a quadratic form in velocities, that the equation locally governing $\Phi(x)$ decouples from the equation for $x(s)$. It follows from the proof of \fullref{Proposition}{Prop1} that this remains true only when $\dot{x}(s)\neq 0$.

\subsection{Solutions of $u''\left(x\right)+h\left(x\right)u\left(x\right)=0$ from geodesic curves on $\mathbb{M}_h$} \label{subsecSol} We are ready to establish the major result of this paper. 
\begin{thm}\label{thm-main2}
Solutions of all linear second order ordinary differential equations (\ref{Eq.main}), with any $h(x)$ of class $C^1$, can be locally expressed through geodesic curves in hyperbolic two dimensional geometry. 
    
More precisely, let  $h(x)$ be of class $C^1$, let $x_0\in\mathbb{R}$ and let 
$\Phi(x)$ be of class $C^2$ on some interval $\left]x_-,x_+\right[$, such that $x_-<x_0<x_+$ and $h\!\left(x\right)\neq\Phi^{2}\!\left(x\right)$ on this interval. If
\begin{subequations}\label{specsol}
\begin{equation}\label{specsol1}
        u_\mathtt{top}(x)=\exp\left[\int_{x_0}^{x}\!\!d\xi\,\Phi\left(\xi\right)\frac{\Phi'\left(\xi\right)-\sqrt{\left[h\left(\xi\right)-\Phi^{2}\left(\xi\right)\right]^{2}+\left[\Phi'\left(\xi\right)\right]^{2}}}{h\left(\xi\right)-\Phi^{2}\left(\xi\right)}\right],
    \end{equation}
    \begin{equation}\label{specsol2}
        u_\mathtt{bot}(x)=\exp\left[\int_{x_0}^{x}\!\!d\xi\,\Phi\left(\xi\right)\frac{\Phi'\left(\xi\right)+\sqrt{\left[h\left(\xi\right)-\Phi^{2}\left(\xi\right)\right]^{2}+\left[\Phi'\left(\xi\right)\right]^{2}}}{h\left(\xi\right)-\Phi^{2}\left(\xi\right)}\right],
    \end{equation}
    \end{subequations}
then for two arbitrary constants $A$ and $B$
\begin{equation}\label{casesS0}
    u(x)=A u_{\mathtt{top}}(x) + B u_{\mathtt{bot}}(x),
\end{equation}
is a general solution of (\ref{Eq.main}) on $\left]x_-,x_+\right[$ if  and only if $\Phi(x)$ is a solution of (\ref{ELexp}) on this interval.    
\end{thm}
We observe that if $\Phi^2(x)$ lays on top/bottom of $h(x)$, the function $u_{\mathtt{top}}(x)$ is increasing/decreasing. Based on \fullref{Theorem}{thm-main2} we immediately get a corollary.
\begin{cor}\label{Cor2} \begin{subequations}\label{qmp}
    Functions
    \begin{equation}\label{qmp1}
        \Theta_{\mathtt{top}}(x)=\frac{u'_{\mathtt{top}}(x)}{u_{\mathtt{top}}(x)}=\Phi\left(x\right)\frac{\Phi'\left(x\right)-\sqrt{\left[h\left(x\right)-\Phi^{2}\left(x\right)\right]^{2}+\left[\Phi'\left(x\right)\right]^{2}}}{h\left(x\right)-\Phi^{2}\left(x\right)},
    \end{equation}
      \begin{equation}\label{qmp2}
        \Theta_{\mathtt{bot}}(x)=\frac{u'_{\mathtt{bot}}(x)}{u_{\mathtt{bot}}(x)}=\Phi\left(x\right)\frac{\Phi'\left(x\right)+\sqrt{\left[h\left(x\right)-\Phi^{2}\left(x\right)\right]^{2}+\left[\Phi'\left(x\right)\right]^{2}}}{h\left(x\right)-\Phi^{2}\left(x\right)},
    \end{equation}\end{subequations}
    are two distinct solutions of a Ricatti equation 
    \begin{equation}\label{Ricatti}
       \Theta'\left(x\right)+\Theta^2\left(x\right)+h(x)=0,
    \end{equation}
if  and only if $\Phi(x)$ is a solution of (\ref{ELexp}) on $\left]x_-,x_+\right[$. Moreover, the relations (\ref{specsol}) and (\ref{qmp}) can be uniquely inverted
 \begin{equation}\label{PhiSqrt}
   \Phi(x) = \sqrt{-\frac{u'_\mathtt{top}(x)u'_\mathtt{bot}(x)}{u_\mathtt{top}(x)u_\mathtt{bot}(x)}} =\sqrt{-\Theta_{\mathtt{top}}(x)\Theta_{\mathtt{bot}}(x)}. 
\end{equation}
\end{cor}
Since $\Phi(x)>0$, the last formula can also be rewritten to the form
 \begin{equation}\label{PhiSqrtGEO}
   \Phi^2(x)u_\mathtt{top}(x)u_\mathtt{bot}(x) +u'_\mathtt{top}(x)u'_\mathtt{bot}(x)=0. 
\end{equation}
Eq. (\ref{PhiSqrt}) turns out to have a more universal character, while looked at from a perspective of a rather obvious statement.
\begin{afa}\label{fact}
    For every $x_0\in\mathbb{R}$, using the freedom of choice of initial conditions, one can specify two linearly independent solutions of (\ref{Eq.main}), such that $u_1\left(x_0\right)\neq 0$, $u_2\left(x_0\right)\neq 0$ and $u'_1\left(x_0\right)u'_2\left(x_0\right)/\left[u_1\left(x_0\right)u_2\left(x_0\right)\right] < 0$. Then, there exists $\left]x_-,x_+\right[$ such that $x_0\in\left]x_-,x_+\right[$ and the solutions $u_1(x)$ and $u_2(x)$ fulfill these conditions on $\left]x_-,x_+\right[$.
\end{afa}
We observe that $u_{\mathtt{top}}(x)$ and $u_{\mathtt{bot}}(x)$ provide an example of such a construction. Moreover, we can prove uniqueness of this choice with respect to $\Phi(x)$, up to multiplicative constants which specify initial conditions at $x_0$.
\begin{prop}\label{prop22}
Let $x_0\in\mathbb{R}$, and let $u_1(x)$ and $u_2(x)$ be a pair  which on $\left]x_-,x_+\right[$ meets the requirements collected in \fullref{Fact}{fact}. Then 
\begin{equation}\label{PhiSqrt12}
   \Phi(x) = \sqrt{-\frac{u'_1(x)u'_2(x)}{u_1(x)u_2(x)}}, 
\end{equation}
is a solution of Eq. (\ref{ELexp}) and 
\begin{equation}\label{casesS}
    u_1(x)=u_1\left(x_0\right) u_{\mathtt{top}/\mathtt{bot}}(x),\qquad u_2(x)=u_2\left(x_0\right) u_{\mathtt{bot}/\mathtt{top}}(x).
\end{equation}
\end{prop}

The opening assertion of \fullref{Theorem}{thm-main2}, even though stated in a slightly bold fashion, refers to the fact that all two dimensional manifolds $\mathbb{M}_h$ are models of the same hyperbolic geometry (see \fullref{Lemma}{lem-hyp}), just expressed in different coordinates specific to $h(x)$ (see \fullref{Lemma}{thm-main} below). In other words, while the explicit form of the solutions of Eq. (\ref{Eq.main}) can tremendously depend on the choice of $h(x)$, its geometric interpretation in terms of geodesic curves is independent of this choice. 

Let us stress that, for any given $h(x)$, special solutions of both (\ref{Eq.main}) and (\ref{ELexp}), provided that both exist, can always be functionally related with each other (at least locally). With a different choice of $h(x)$ such a relation is generally expected to be of a completely different form. The same is true for a single choice of $h(x)$ and distinct parts of the domain. However, as \fullref{Theorem}{thm-main2} shows, the relations (\ref{specsol}) have a universal character. This is a very peculiar and unexpected result, which uncovers a fundamental connection between  linear second order ordinary differential equations and hyperbolic geometry in two dimensions.

It is also clear that the formulas (\ref{specsol}) cannot in a general case provide global solutions, i.e. we do not get $x_\pm=\pm\infty$. For example, if $h(x)$ is a positive constant, it is well known that every solution $u(x)$ reaches negative values. On the other hand, all solutions in (\ref{specsol}) are positive. In fact, the same is true while reducing (\ref{Eq.main}) to the Ricatti equation (\ref{Ricatti}) with the help of the substitution $u'\left(x\right)=\Theta\left(x\right)u\left(x\right)$. We examine that feature while describing intuitive examples in \fullref{Sec.}{subsecEX}. 

\subsection{Diffeomorphism between $\mathbb{M}_h$ and Poincar{\'e} upper half plane} After the major result we move back to geometric aspects of the Riemannian manifold $\mathbb{M}_h$, previously discussed in \fullref{Sec.}{subsecCur} and \fullref{Sec.}{subsecGeo}. 
\begin{lem}
    \label{thm-main}
For any pair $u_1(x)$ and $u_2(x)$ of real-valued linearly independent solutions of (\ref{Eq.main}) with $h(x)$ of class $C^0$, the map $\mathcal{H}\rightarrow\mathcal{H}$:
\begin{subequations}\label{map}
\begin{equation}\label{DiffX}
   X\left(x,\varPhi\right)= X_0\pm W^{-1}\frac{\varPhi^2\, u_{1}\left(x\right)u_{2}\left(x\right)+u'_{1}\left(x\right)u'_{2}\left(x\right)}{\varPhi^2\left[u_{1}\left(x\right)\right]^{2}+\left[u'_{1}\left(x\right)\right]^{2}}, 
\end{equation}
\begin{equation}\label{DiffY}
    Y\left(x,\varPhi\right)=\frac{\varPhi}{\varPhi^2\left[u_{1}\left(x\right)\right]^{2}+\left[u'_{1}\left(x\right)\right]^{2}},
\end{equation}
\end{subequations}
is a local diffeomorphism $\mathbb{M}_h\rightarrow\mathbb{H}$. Both $X_0\in\mathbb{R}$ and the $\pm$ sign in (\ref{DiffX}) are arbitrary. The Wronskian $W= u'_{1}(x)u_{2}(x)-u_{1}(x)u'_{2}(x)$ is a normalization constant which, by assumption of linear independence, is different than zero. 
\end{lem}
The Jacobian of this coordinate transformation vanishes only when  $\varPhi^2=h(x)$, so only outside of $\mathcal{H}_0$. Let us stress that since $\varPhi>0$, while $u_{1}(x)$ and $u'_{1}(x)$ cannot vanish simultaneously, the denominator in (\ref{DiffX}) and (\ref{DiffY}) is never zero. Looking at (\ref{DiffY}) we observe that only the solution $u_1(x)$ has been used to define $Y\left(x,\varPhi\right)$. Though, this is just a mere choice of convention. In fact, since the result works for an arbitrary pair of linearly independent solutions, we get:
\begin{cor}\label{Cor1}
    Group $GL(2,\mathbb{R})$ acting on the vector space spanned by linearly independent solutions of (\ref{Eq.main}) induces a local isomorphism of $\mathbb{H}$. Moreover, $Y\left(x,\varPhi\right)$ defined in (\ref{DiffY}) with $u_1(x)$ replaced by a general solution of (\ref{Eq.main}), in the form $u(x)=A u_1(x)+B u_2(x)$ with arbitrary real constants $A$ and $B$, is a complete integral of a nonlinear and non-separable partial differential equation
\begin{equation}\label{pardif}
\left(\frac{\varPhi}{h\left(x\right)-\varPhi^{2}}\frac{\partial}{\partial x}Y\left(x,\varPhi\right)\right)^{2}+\left(\varPhi\frac{\partial}{\partial\varPhi}Y\left(x,\varPhi\right)\right)^{2}=Y^2\left(x,\varPhi\right),
\end{equation}
involving an arbitrary continuous function $h(x)$.     
\end{cor}
The requirement of non-vanishing determinant, defining  the $GL\left(2,\mathbb{R}\right)$ group, assures linear independence of all new linear combinations of the solutions. Partial differential equation (\ref{pardif}) means that a gradient of $\ln Y\left(x,\varPhi\right)$ has a unit norm in $\mathbb{M}_h$, the same way as it is trivially true in $\mathbb{H}$. Even though solving (\ref{pardif}) has only been a step on a way to arrive at \fullref{Lemma}{thm-main},  this solution by itself might turn out to be useful elsewhere.

It is well known that geodesic curves on $\mathbb{H}$ are either vertical straight lines, or semi circles. From (\ref{PhiSqrtGEO}) and (\ref{DiffX}) we can immediately observe that if we pick a geodesic curve $\Phi(x)$ and set $u_1(x)=u_{\mathtt{top}}(x)$ and $u_2(x)=u_{\mathtt{bot}}(x)$ (or the other way around), then this curve is mapped to $X=X_0$. This is a vertical line. In fact, all geodesic curves on $\mathbb{H}$ can be solved to give explicit relations between $\Phi(x)$, $u_1(x)$, and $u_2(x)$. Even though, for the current purpose, such results do not bring more insight in comparison with the above vertical solution related with the formulas (\ref{specsol}), I sketch a less tedious way leading in that direction. To this end, we recall that $\mathbb{H}$ possesses three Killing vectors (the maximal set): $k_1 = \left(X^2-Y^2\right)\partial_X+2 X Y \partial_Y$, $k_2=X \partial_X+Y \partial_Y$ and $k_3=\partial_X$, so that there exist three constants of motion $A_i$, $i=1,2,3$, associated with each Killing vector $k_i$. With the help of \fullref{Lemma}{thm-main}, after finding the Jacobian matrix $J$ of (\ref{map}), we can immediately get three Killing vectors in $\mathbb{M}_h$, which are $J^{-1}k_i$. It is then straightforward to write down three conserved quantities associated with these Killing vectors and equal to $A_i$. One can algebraically solve any two of them, getting $\dot x(s)$ and $\dot\varPhi(s)$ as functions of $x(s)$ and $\varPhi(s)$. The third constant of motion then leads to the relation
\begin{equation}
\Phi^2\left(x\right)+\frac{A_{1}W^{2}\left[u_{1}'\left(x\right)\right]^{2}-2A_{2}W u_{1}'\left(x\right)u_{2}'\left(x\right)+A_{3}\left[u_{2}'\left(x\right)\right]^{2}}{A_{1}W^{2}u_{1}^{2}\left(x\right)-2A_{2}W u_{1}\left(x\right)u_{2}\left(x\right)+A_{3}u_{2}^{2}\left(x\right)}=0.
\end{equation}
Our major result is associated with the special choice $A_1=0=A_3$ and $A_2\neq0$.

\subsection{Examples} \label{subsecEX}
At the beginning of the chapter about linear second order ODEs, Hartman lists three simplest cases \cite{Hartman}: $h(x)=0$, $h(x)=-\omega^2$ and $h(x)=+\omega^2$, with $\omega\neq 0$ being a constant. Without loss of generality we may take $\omega>0$. In order to illustrate \fullref{Theorem}{thm-main2} and \fullref{Lemma}{thm-main} let us examine solutions for these three cases in detail. Special attention will be put on the domain where the solutions are positive. Even in such a basic setting we will be able to observe interesting aspects of the developed formalism. For simplicity we set $X_0=0$ in (\ref{DiffX}). Details of the examples can be found in \fullref{Sec.}{Sec5}.
\begin{notation}[Labeling inside examples]
The Riemannian manifold $\mathbb{M}_h$, for specific choices of $h(x)$ listed above, will be denoted as: $\mathbb{M}_0$ for $h(x)=0$, $\mathbb{M}_-$ for $h(x)=-\omega^2$ and $\mathbb{M}_+$ for $h(x)=+\omega^2$. In the same vein, we denote: $\Phi_{0}(x)$, $u_{1,0}(x)$, $u_{2,0}(x)$, $W_0$, $ u_{\mathtt{top}/\mathtt{bot},0}(x)$, $\Theta_{\mathtt{top}/\mathtt{bot},0}(x)$, $W_{\mathtt{t-b},0}$, etc.
\end{notation}

\begin{exa}\label{exa1}
For $h(x)=0$, the general nonnegative solution of (\ref{ELexp}) is
\begin{equation}\label{Phi0ex1}
    \Phi_{0}\left(x\right)=\frac{1}{\sqrt{C_{1}^2-\left(x+C_{2}\right)^{2}}},
\end{equation}
where $C_1,C_2$ are arbitrary real integration constants and $\left|x+C_{2}\right|\leq\left|C_{1}\right|$. From this inequality we can reconstruct the values $x_\pm=\pm\left|C_{1}\right|-C_2$. Varying both integration constants we are able to establish solutions in any desired interval on a real line. 

Consistently with standard results, the solutions (\ref{specsol}) are linear functions 
\begin{equation}\label{specsollinear}
        u_{\mathtt{top},0}(x)=\frac{x-x_{-}}{x_{0}-x_{-}},\quad u_{\mathtt{bot},0}(x)=\frac{x-x_{+}}{x_{0}-x_{+}},\quad W_{\mathtt{t-b},0}=\frac{x_+-x_-}{\left(x_+-x_0\right)\left(x_0-x_-\right)},
\end{equation}
However, knowing that both $x,x_0\in\left]x_-,x_+\right[$, we observe that numerator and denominator in $u_{\mathtt{top},0}(x)$ are both positive, while in  $u_{\mathtt{bot},0}(x)$ they are both negative. Consequently, both solutions are positive and linearly independent since they represent lines with positive and negative "velocity" respectively. The latter also means that $u'_{\mathtt{top},0}(x)u'_{\mathtt{bot},0}(x)<0$ on $\left]x_-,x_+\right[$.

The above solutions can be used to establish a local diffeomorphism between $\mathbb{M}_0$ and $\mathbb{H}$, however, for finite values of $C_1$ and $C_2$ it will never be global. Still, the manifold $\mathbb{M}_0$ is so simple that we recognize there is a basic global diffeomorphism $X\left(x,\varPhi\right)=x$, $Y\left(x,\varPhi\right)=\varPhi^{-1}$. It can trivially be realized by \fullref{Lemma}{thm-main}, letting $u_{1,0}(x)=1$ and $u_{2,0}(x)=x$, for which $W_0 = -1$, and taking the minus sign in (\ref{DiffX}).
\end{exa}

\begin{exa}\label{exa2}
For $h(x)=-\omega^2$  we can check that $\Phi_{-}(x)=\omega$ is a solution of (\ref{ELexp}) which exists for all $x\in\mathbb{R}$.  From \fullref{Theorem}{thm-main2}, consistently with standard results, we then find
\begin{equation}\label{solminus}
        u_{\mathtt{top},-}(x)=e^{\omega x},\quad u_{\mathtt{bot},-}(x)=e^{-\omega x},
\end{equation}
where for simplicity we took $x_0=0$. These solutions are global.

We  observe that the  metric (\ref{metric}) is nondegenerate on $\mathcal{H}$, so that $\mathcal{H}_0=\mathcal{H}$. In \fullref{Lemma}{thm-main} we can use (\ref{solminus}), for which $W_{\mathtt{t-b},-}=2\omega$,
getting the map
\begin{subequations}\label{glodiff}
\begin{equation}\label{regularexp}
X\left(x,\varPhi\right)=\frac{e^{-2\omega x}}{2\omega}\frac{\varPhi^2-\omega^2}{\varPhi^2+\omega^2},\quad Y\left(x,\varPhi\right)=\varPhi\frac{ e^{-2\omega x}}{\varPhi^2+\omega^{2}}.
\end{equation}
We took plus sign in (\ref{DiffX}). This map is bijective in the whole domain $\mathcal{H}$, and remains such even if extended to a map $\mathbb{R}^2\rightarrow\mathbb{R}^2$. Its global inverse is
\begin{equation}\label{inverseexp}
    x\left(X,Y\right)=-\frac{\ln\left(2\omega\sqrt{X^2+Y^2}\right)}{2\omega},\quad \varPhi\left(X,Y\right)=\omega\frac{X+\sqrt{X^{2}+Y^{2}}}{Y}.
\end{equation}
\end{subequations}
We can see that (\ref{glodiff}) is a global diffeomorphism between $\mathbb{H}$ and $\mathbb{M}_{-}$. Note that, even in this rather simple case, the diffeomorphism is not elementary.
\end{exa}

\begin{exa}\label{HOexa} In the former examples we started from a solution of the geodesic equation (\ref{ELexp}) in order to construct special linearly independent solutions of (\ref{Eq.main}), which are given by (\ref{specsol1}) and (\ref{specsol2}). Now, we start with two preselected solutions of (\ref{Eq.main}) and will construct the solution (\ref{PhiSqrt12}) of the geodesic equation which, as also showed in \fullref{Proposition}{prop22}, is the same as  (\ref{PhiSqrt}).

We consider harmonic oscillator given by $h(x)=+\omega^2$. This model is somehow special, because standard solutions
\begin{equation}\label{HOB}
        u_{1,+}(x)=\cos\left[\omega\left(x-\bar{x}\right)\right],\quad u_{2,+}(x)=\sin\left[\omega\left(x-\bar{x}\right)\right],\quad W_+=-\omega,
\end{equation}
with an arbitrary shift $\bar{x}$, even though for $k\in \mathbb{Z}$ and $\bar{x}+k\pi/(2\omega)<x_0(k)<\bar{x}+(k+1)\pi/(2\omega)$ meet the requirements listed in \fullref{Fact}{fact}, give $\Phi_+(x)=\omega$ while plugged into (\ref{PhiSqrt12}). This is obviously singular (does not belong to $\mathcal{H}_0$). A resolution of that problem is to consider the solutions (\ref{HOB}) with a relative phase shift $0<\tilde{x}<\pi/(2\omega)$:
\begin{equation}
        u_{1,+}(x),\qquad u_{+,2}(x-\tilde{x}),\quad W_{+}\left(\tilde{x}\right)=-\omega\cos\left(\omega\tilde{x}\right).
\end{equation}
From (\ref{PhiSqrt12}) we then get the solution of (\ref{ELexp})
\begin{equation}\label{solcot}
\Phi_+\left(x\right)=\omega\sqrt{\cot\left[\omega\left(x-\bar{x}-\tilde{x}\right)\right]\tan\left[\omega\left(x-\bar{x}\right)\right]},
\end{equation}
while due to periodicity of trigonometric functions we read out infinitely many domains of existence $\left]x_-(k),x_+(k)\right[$ with 
\begin{equation}\label{domainwithk}
    x_-(k)=\frac{k \pi}{2\omega} + \bar{x}+\tilde{x},\qquad x_+(k)=\bar{x}+\frac{(k+1) \pi}{2\omega},\quad k\in \mathbb{Z}.
\end{equation}
We observe that $x_-(k)$ increased by $\tilde{x}$ in comparison with the domain for (\ref{HOB}).

Interestingly, the singularity $\Phi_{+}(x)=\omega$ never occurs because $\Phi_+\left(x\right)>\omega$ in the selected domain for even values of $k$, while for odd values of $k$ we get the case $0<\Phi_+\left(x\right)<\omega$. This boils down to one of the cases in (\ref{casesS}).

Even though basic, this example is topologically non-trivial since $\mathbb{M}_{+}$ is a disjoint union of a strip $\{\left(x,\varPhi\right)\in\mathbb{R}^2\,|\,0<\varPhi<\omega\}$ and a half plane $\{\left(x,\varPhi\right)\in\mathbb{R}^2\,|\,\varPhi>\omega\}$. In \fullref{Lemma}{thm-main} we will use (\ref{HOB}), getting the transformation
\begin{equation}
X\left(x,\varPhi\right)=\frac{\left(\omega^{2}-\varPhi^{2}\right)\sqrt{\zeta(x)\left(1-\zeta(x)\right)}}{\omega\left[\varPhi^{2}\zeta(x)+\omega^{2}\left(1-\zeta(x)\right)\right]},\quad Y\left(x,\varPhi\right)=\frac{\varPhi}{\varPhi^{2}\zeta(x)+\omega^{2}\left(1-\zeta(x)\right)},\label{XHO}
\end{equation}
where $\zeta\left(x\right)=\cos^{2}\left[\omega\left(x-\bar{x}\right)\right]$. Again we took plus sign in (\ref{DiffX}).

If $x\in\left]\bar{x}+k\pi/(2\omega),\bar{x}+(k+1)\pi/\left(2\omega\right)\right[$, for $k\in \mathbb{Z}$, injectivity of that map will be assured. Without loss of generality we set $k=0$. Despite a more complicated form than in the previous example, this map is tractable because for fixed $\varPhi<\omega$ its image turns out to be (see \fullref{Sec.}{Sec5} for details) a right semicircle ($X>0$) of
\begin{equation}\label{circleHOM}
X^{2}+\left(\frac{\varPhi^{2}+\omega^{2}}{2\varPhi\omega^{2}}-Y\right)^{2}=\left(\frac{\varPhi^{2}-\omega^{2}}{2\varPhi\omega^{2}}\right)^{2},
\end{equation}
while the left part $X<0$ furnishes  $\varPhi>\omega$. Moreover, again for fixed $\varPhi$, the function $Y\left(x,\varPhi\right)$ from (\ref{XHO}) is monotonic in the selected open-interval domain of $x$ and with its values covers the entire semicircle corresponding to $\varPhi$. In other words, for any $\bar{x}\in\mathbb{R}$ we get a single diffeomorphism
\begin{equation}\label{domM}
\mathbb{M}_{+}\cap\mathcal{R}^{<}_{\bar{x}}\rightarrow\mathbb{H}\cap\mathcal{Q}_{1},\qquad\mathbb{M}_{+}\cap\mathcal{R}^{>}_{\bar{x}}\rightarrow\mathbb{H}\cap\mathcal{Q}_{2},
\end{equation}
where $\mathcal{R}^{<}_{\bar{x}}=\left]\bar{x},\bar{x}+\pi/\left(2\omega\right)\right[\times\left]0,\omega\right[$
and $\mathcal{R}^{>}_{\bar{x}}=\left]\bar{x},\bar{x}+\pi/\left(2\omega\right)\right[\times\left]\omega,\infty\right[$,
and $\mathcal{Q}_{k}$, for $k=1,\ldots,4$ are standard  quadrants of the
plane $\mathbb{R}^{2}$ taken without the ordinate (for example $\mathcal{Q}_{1}=\left]0,\infty\right[\times\mathbb{R}$). We note that all involved sets: $\mathcal{H}\cap\mathcal{R}^{<}_{\bar{x}}$, $\mathcal{H}\cap\mathcal{R}^{>}_{\bar{x}}$, $\mathcal{H}\cap\mathcal{Q}_{1}$ and $\mathcal{H}\cap\mathcal{Q}_{2}$ are open.
\end{exa}

\section{Discussion}
\label{Sec3}
Let me summarize the main message of this contribution. While it would not be astonishing that the second order ODE (\ref{Eq.main}) with a fixed choice of $h(x)$ is equivalent to a geodesic equation for some specific metric, it shall be considered remarkable that such a metric can be chosen in a universal way (\ref{metric}), being a basic function of $h(x)$. Even more remarkably, for any $h(x)$ the metric (\ref{metric}) describes the same profound case of hyperbolic geometry,  modulo topological peculiarities. Conjunction of these two properties gives a meaning to the presented results.

While pondering on the above statement, as an exercise,  consider a flat geometry with metric $dv^2+du^2$. It admits a unit-speed geodesic line (in explicit form) $u\left(v\right)=v$. By a coordinate transformation $v=v\left(x\right)$, with an arbitrary invertible function $v(x)$, the metric transforms to $\left[V\left(x\right)\right]^2dx^2+du^2$, where $V(x)=v'(x)$,
while the geodesic line becomes $u\left(x\right)=v\left(x\right)$. In this way, one can convince themselves that freedom of coordinate transformations might be enough to let the geodesic curve become an arbitrary function, at least locally. We can also see that the geodesic curve in question fulfills the most basic differential equation $u'(x)=V(x)$. Solutions to this equation, which using a substitution $u(x)=\ln{\tilde{u}(x)}$ becomes equivalent to a linear homogeneous first order ODE $\tilde{u}'(x)=V(x)\tilde{u}(x)$, locally are equal to geodesic curves in flat space. However, if we wish to go one step further, and let this in fact  arbitrary function $u(x)$ obey (\ref{Eq.main}), we simply get $v''(x)+h(x)v(x)=0$. So, to use the above approach to propose a flat two dimensional metric for which its geodesic curve is a solution of (\ref{Eq.main}), given an input function $h(x)$, we would first need to solve (\ref{Eq.main}) itself. As a consequence, $V(x)$ cannot in general be an explicit function (or even functional) of $h(x)$. Surprisingly, if flat geometry is replaced by hyperbolic geometry, an analogous result holds. 

In other words, the above flat geometry and solutions of the differential equation (\ref{Eq.main}) are not \textit{independent} entities which "unify to our surprise", since both \textit{a priori} require solving (\ref{Eq.main}). On the other hand, in our hyperbolic scenario, geometry of the Riemannian manifold $\mathbb{M}_h$ has independence. One can study its properties, for example a number of disjoint pieces, without at all referring to (\ref{Eq.main}). Only \textit{a posteriori} we discover a deep mutual connection between $\mathbb{M}_h$ and Eq.  (\ref{Eq.main}). Therefore, the results derived in this paper contribute to the very center of the theory of linear second order ODEs, a field, which otherwise is customarily assumed to be very well-explored. In a way, we supplement the theory of ODEs at a level perhaps more suited for the turn of the 20st century. In that regard, it is fair to say that the original motivation behind the current study comes from theory of Shabat-Zakharov systems \cite{ShabatZakharov}. With some adjustments, I have initially reinterpreted an upper bound on the reflection coefficient \cite{Visser} in geometric terms, recognizing it is a function of an arc length of a (non-geodesic) curve on a plane equipped with metric $g_h$. Without this observation I would neither define $\mathbb{M}_h$ nor explore its geometry looking for strict relations with solutions of (\ref{Eq.main}).

To recapitulate the technical findings of this manuscript, equivalence between real nonsingular linear second order ordinary differential equations and a geodesic equation in two dimensional hyperbolic geometry shall be understood by means of two major results proved here. \fullref{Theorem}{thm-main2} shows that a solution of the geodesic equation in $\mathbb{M}_h$ determines two linearly independent solutions of (\ref{Eq.main}), while due to \fullref{Proposition}{prop22} any two linearly independent solutions of (\ref{Eq.main}) determine a solution of the geodesic equation on an interval in which both aforementioned solutions of (\ref{Eq.main}) are positive and growing in opposite directions. \fullref{Lemma}{thm-main} shows that the Riemannian manifold $\mathbb{M}_h$ is locally diffeomorphic to the Poincar{\'e} upper half plane. 

Let me also list a few directions of potential technical extensions. First of all, the question remains about types of singularities, discontinuities, zeros, non-differentiability, etc. in the original equation (\ref{Eq.1}), which allow to keep the geometric link established in this paper. Clearly, if there is a subdomain on the real line in which none of these issues occur, the formalism works locally. Also, the case of removable singularities shall not pose a fundamental problem. 

Being more precise, of relevance is a careful study of minimal required continuity assumptions about the function $h(x)$. Formally speaking, for $h(x)$ being of class $C^2$ everything works well, however, geometric properties of the metric (\ref{metric}) are intact even if $h(x)$ is only absolutely continuous. On the other hand, to consider the geodesic equation we generically need $h(x)$ to be of class $C^1$. Moreover, what happens when required continuity property of $h(x)$ only fails at a point $x$ in which $u(x)=0$?

Through the aforementioned arc length functional there is a natural connection with an old inverse problem of the calculus of variations. Even though its solution for Eq. (\ref{Eq.main}) is obvious, the current approach offers an indirect alternative. In that regard it is worth a consideration if other non-linear second (or higher) order ODEs (e.g. those stemming from variational principles) can also be linked with particular types of geometry, and classified according to properties such as curvature.

Next, despite a lack of a hard evidence, I anticipate a deeper link between the current findings and a notion of the Schwarzian derivative. The latter one is strongly linked with M{\"o}bius transformations acting on a hyperbolic plane. 

For the end of the Discussion part I leave a natural big open question, which pertains to the complex scenario. If the real variable $x$ is replaced by its complex counterpart $z$, and all involved functions are allowed to be complex, the real dimensionality of the problem doubles. This suggests that the hyperbolic pattern shall not in general hold for complex variants of (\ref{Eq.main}). However, we are still able to give a distinct geometric counterpart of the differential equation under consideration. The sequel paper is fully devoted to this issue.

\section{Proofs}
\label{Sec4}
Even though proofs presented below are of a rather computational type, they were nailed down to the form in which they can be reproduced with pen and paper. However, in the case of literally every obtained result, the road to find settings which in an easy way show the merit, was far more involved (for example, see \fullref{Corollary}{Cor1} and discussion below). A few more results of a similar level of technical difficulty were obtained, but are not going to be reported here since they are not essential to understand the current paper. From now on we start one-by-one proofs of particular statements.

\begin{proof}[Proof of \textup{\fullref{Lemma}{lem-hyp}}]
We are going to describe the geometry of $g_{h}$ by providing all
non-vanishing components of associated Christoffel symbols and Riemann tensor, both defined in a standard textbook way. Since conventions might differ, for consistency
I provide explicit expressions for these objects in some coordinates (Einstein summation convention applies) \cite{Besse}: 
\begin{subequations}
\begin{equation}
\Gamma_{\:\:\:jk}^{i}=\frac{1}{2}g^{il}\left(g_{lj,k}+g_{lk,j}-g_{jk,l}\right),
\end{equation}
\begin{equation}
R_{\:\:\:jkl}^{i}=\Gamma_{\:\:\:jl,k}^{i}-\Gamma_{\:\:\:jk,l}^{i}+\Gamma_{\:\:\:mk}^{i}\Gamma_{\:\:\:\:jl}^{m}-\Gamma_{\:\:\:ml}^{i}\Gamma_{\:\:\:\:jk}^{m}.
\end{equation}
\end{subequations}
In the chosen
coordinates $\left(x,\varPhi\right)$ the metric $g_h$ is (I skip the label $h$)
\begin{subequations}
\begin{equation}
    g_{xx} = \frac{\left(h\left(x\right)-\varPhi^{2}\right)^{2}}{\varPhi^{2}},\quad g_{\varPhi\varPhi}=\frac{1}{\varPhi^{2}},\quad g_{x\varPhi}=g_{\varPhi x}=0.
\end{equation}
Therefore:
\begin{equation}
\Gamma_{\:\:\:xx}^{x}=\frac{h'\left(x\right)}{h\left(x\right)-\varPhi^{2}},\quad\Gamma_{\:\:\:x\varPhi}^{x}=\Gamma_{\:\:\:\varPhi x}^{x}=\frac{\varPhi^{2}+h\left(x\right)}{\varPhi\left[\varPhi^{2}-h\left(x\right)\right]},
\end{equation}
\begin{equation}
\Gamma_{\:\:\:xx}^{\varPhi}=\frac{h^{2}\left(x\right)-\varPhi^{4}}{\varPhi},\quad\Gamma_{\:\:\:\varPhi\varPhi}^{\varPhi}=-\frac{1}{\varPhi},
\end{equation}
\begin{equation}
R_{\:\:\:\varPhi x\varPhi}^{x}=-R_{\:\:\:\varPhi\varPhi x}^{x}=-\frac{1}{\varPhi^{2}},\quad R_{\:\:\:xx\varPhi}^{\varPhi}=-R_{\:\:\:x\varPhi x}^{\varPhi}=\frac{\left(h\left(x\right)-\varPhi^{2}\right)^{2}}{\varPhi^{2}},
\end{equation}
\begin{equation}
R_{x\varPhi x\varPhi}=-R_{x\varPhi\varPhi x}=-R_{\varPhi xx\varPhi}=R_{\varPhi x\varPhi x}=-\frac{\left(h\left(x\right)-\varPhi^{2}\right)^{2}}{\varPhi^{4}}.
\end{equation}
I have only listed non-vanishing coefficients. The Ricci tensor $R_{ij}=R_{\:\:\:ikj}^{k}$ is
\begin{equation}
R_{xx}=-\frac{\left(h\left(x\right)-\varPhi^{2}\right)^{2}}{\varPhi^{2}},\quad R_{\varPhi\varPhi}=-\frac{1}{\varPhi^{2}},\quad R_{x\varPhi}=R_{\varPhi x}=0.
\end{equation}    
\end{subequations}
From the last equation we notice that the Ricci tensor is $-g_{h}$.
Since $\textrm{Tr}g_{h}=2$, we get the value of the Ricci scalar equal
to $-2$. In two dimensions this already sets the value of $-1$ for
sectional curvature, equal to the Gaussian curvature. Obviously, the
same follows from the formula
\begin{equation}
K=\frac{R_{x\varPhi x\varPhi}}{g_{xx}g_{\varPhi\varPhi}-g_{x\varPhi}g_{\varPhi x}}\equiv-1.
\end{equation}

Concerning the accompanying remark, we can see that $R_{x\varPhi x\varPhi}$
does only depend on $h\left(x\right)$ and terms involving its derivatives
are removed due to antisymmetry. In addition, the Riemann tensor is
regular at $\varPhi^{2}=h\left(x\right)$.
\end{proof}

\begin{proof}[Proof of \textup{\fullref{Proposition}{Prop1}}]
We start from the substitution $\varPhi\left(s\right)=\Phi\left(x\left(s\right)\right)$,
which leads to:\begin{subequations}
\begin{equation}\label{primeph}
\dot{\varPhi}\left(s\right)=\Phi'\left(x\left(s\right)\right)\dot{x}\left(s\right),
\end{equation}
\begin{equation}
\ddot{\varPhi}\left(s\right)=\Phi'\left(x\left(s\right)\right)\ddot{x}\left(s\right)+\Phi''\left(x\left(s\right)\right)\left[\dot{x}\left(s\right)\right]^{2}.
\end{equation}\end{subequations}
From now on we drop the arguments $s$ and $x(s)$, recalling that $\Phi$ is a function of $x$ while $\varPhi$ is a function of $s$, as well as that the prime denotes the derivative with respect to $x$ while the dot is the derivative with respect to $s$. We substitute $\ddot x$ and $\ddot \varPhi$ from the geodesic equations (\ref{geode}) getting
\begin{equation}
\frac{\Phi^{4}-h^{2}}{\Phi}\dot{x}^{2}+\frac{\dot{\varPhi}^{2}}{\Phi}=\Phi'\left(\frac{h'}{\Phi^{2}-h}\dot{x}^{2}-2\frac{\Phi^{2}-h}{\Phi\left(\Phi^{2}-h\right)}\dot{x}\dot{\varPhi}\right)+\Phi''\dot{x}^{2}.
\end{equation}
After substituting $\dot{\varPhi}$ from (\ref{primeph}) we collect together two terms with $\left[\Phi'\right]^{2}$. After moving all the terms to the right hand side we get
\begin{equation}\label{finalhelpf}
0=\dot{x}^2\left\{\Phi''-\frac{\Phi^{4}-h^{2}}{\Phi}+\frac{h'\Phi'}{\left(\Phi^{2}-h\right)}-\frac{\left(3\Phi^{2}+h\right)\left[\Phi'\right]^{2}}{\Phi\left(\Phi^{2}-h\right)}\right\}.
\end{equation}
As long as $\dot{x}\neq0$, this is the same equation as Eq. (\ref{ELexp}). By direct inspection we can confirm that, assuming  $\dot{x}\neq0$ and using the Christoffel symbols listed in the proof of \fullref{Lemma}{lem-hyp}, Eq. (\ref{finalhelpf}) can be rewritten to the form
\begin{equation}
\Phi''=-\Gamma_{\:\:\:xx}^{\varPhi}+\Gamma_{\:\:\:xx}^{x}\Phi'+\left(2\Gamma_{\:\:\:\varPhi x}^{x}-\Gamma_{\:\:\:\varPhi\varPhi}^{\varPhi}\right)\left[\Phi'\right]^{2}.
\end{equation}
\end{proof}

\begin{proof}[Proof of \textup{\fullref{Theorem}{thm-main2}}]
The same way as in the former proof  we omit the dependence on the variable $x$. The same applies to proofs of  \fullref{Corollary}{Cor2} and \fullref{Proposition}{prop22}. In order to simplify the derivation we introduce a norm of a velocity of a geodesic curve in explicit form \begin{subequations}
\begin{equation}
\mathcal{L}=\sqrt{\left(h-\Phi^{2}\right)^{2}+\left[\Phi'\right]^{2}},
\end{equation}
and denote
\begin{equation}\label{Thetapm}
\Theta_{\pm}=\Phi\frac{\Phi'\pm\mathcal{L}}{h-\Phi^{2}}.
\end{equation}
We find
\begin{equation}
\Theta_{\pm}'=\Phi'\frac{\Phi'\pm\mathcal{L}}{h-\Phi^{2}}+\Phi\frac{\left(\Phi''\pm\mathcal{L}'\right)\left(h-\Phi^{2}\right)-\left(\Phi'\pm\mathcal{L}\right)\left(h'-2\Phi\Phi'\right)}{\left(h-\Phi^{2}\right)^{2}},
\end{equation}
and
\begin{equation}
\Theta_{\pm}^{2}=2\Phi^{2}\Phi'\frac{\Phi'\pm\mathcal{L}}{\left(h-\Phi^{2}\right)^{2}}+\Phi^{2}.
\end{equation}    
\end{subequations}\begin{subequations}
We then get
\begin{equation}\label{RicPQ}
\Theta_{\pm}'+\Theta_{\pm}^{2}+h=P\pm Q,
\end{equation}
where
\begin{equation}
P=\frac{\left[\Phi'\right]^{2}+\Phi\Phi''}{h-\Phi^{2}}-\Phi\Phi'\frac{h'-4\Phi\Phi'}{\left(h-\Phi^{2}\right)^{2}}+\Phi^{2}+h,
\end{equation}
\begin{equation}
Q=\frac{\Phi'\mathcal{L}+\Phi\mathcal{L}'}{h-\Phi^{2}}-\Phi\frac{h'-4\Phi\Phi'}{\left(h-\Phi^{2}\right)^{2}}\mathcal{L}.
\end{equation}\end{subequations}

We first rearrange\begin{subequations}
\begin{eqnarray}
P & = & \frac{\left[\Phi'\right]^{2}+\Phi\Phi''}{h-\Phi^{2}}-\Phi\Phi'\frac{h'-4\Phi\Phi'}{\left(h-\Phi^{2}\right)^{2}}+\Phi^{2}+h\nonumber \\
 & = & \frac{\Phi\Phi''}{h-\Phi^{2}}+\Phi'\frac{\Phi'\left(h-\Phi^{2}\right)-\Phi\left(h'-4\Phi\Phi'\right)}{\left(h-\Phi^{2}\right)^{2}}+\Phi^{2}+h\nonumber \\
 & = & \frac{\Phi\Phi''}{h-\Phi^{2}}+\Phi'\frac{\Phi'\left(h+3\Phi^{2}\right)-\Phi h'}{\left(h-\Phi^{2}\right)^{2}}+\Phi^{2}+h,\label{Pfinal}
\end{eqnarray}
and
\begin{equation}
Q=\frac{\Phi}{h-\Phi^{2}}\mathcal{L}'+\frac{\Phi'\left(h+3\Phi^{2}\right)-\Phi h'}{\left(h-\Phi^{2}\right)^{2}}\mathcal{L}.
\end{equation}\end{subequations}
Next, since\begin{subequations}
\begin{equation}
\mathcal{L}'=\frac{\left(h-\Phi^{2}\right)\left(h'-2\Phi\Phi'\right)+\Phi'\Phi''}{\mathcal{L}},
\end{equation}
we observe that $Q/\Phi'$ is regular at $\Phi'= 0$. Therefore, we can write
\begin{eqnarray}
\frac{\mathcal{L}Q}{\Phi'} & = & \frac{\Phi}{h-\Phi^{2}}\frac{\mathcal{L}\mathcal{L}'}{\Phi'}+\frac{h\Phi'+\Phi\left(3\Phi\Phi'-h'\right)}{\left(h-\Phi^{2}\right)^{2}}\frac{\mathcal{L}^{2}}{\Phi'}\nonumber \\
 & = & \frac{\Phi\Phi''}{h-\Phi^{2}}+\frac{h'-2\Phi\Phi'}{\Phi^{-1}\Phi'}+\frac{\Phi'\left(h+3\Phi^{2}\right)-\Phi h'}{\left(h-\Phi^{2}\right)^{2}}\left\{ \frac{\left(h-\Phi^{2}\right)^{2}}{\Phi'}+\Phi'\right\} \nonumber \\
 & = & \frac{\Phi\Phi''}{h-\Phi^{2}}+\frac{\bcancel{h'}-\cancel{2\Phi\Phi'}}{\xcancel{\Phi^{-1}\Phi'}}+\left\{ h+\cancel{3}\Phi^{2}-\xcancel{\frac{\Phi}{\Phi'}}\bcancel{h'}+\frac{\Phi'\left(h+3\Phi^{2}\right)-\Phi h'}{\left(h-\Phi^{2}\right)^{2}}\Phi'\right\} \nonumber \\
 & = & \frac{\Phi\Phi''}{h-\Phi^{2}}+h+\Phi^{2}+\frac{\Phi'\left(h+3\Phi^{2}\right)-\Phi h'}{\left(h-\Phi^{2}\right)^{2}}\Phi'\nonumber \\ \label{QP}
 & = & P,
\end{eqnarray}
where we used
\begin{eqnarray}
\frac{\Phi}{h-\Phi^{2}}\frac{\mathcal{L}\mathcal{L}'}{\Phi'}&=&\frac{\Phi\left[\left(h-\Phi^{2}\right)\left(h'-2\Phi\Phi'\right)+\Phi'\Phi''\right]}{\Phi'\left(h-\Phi^{2}\right)}\\ \nonumber
& =&\frac{\Phi\Phi''}{h-\Phi^{2}}+\frac{h'-2\Phi\Phi'}{\Phi^{-1}\Phi'},
\end{eqnarray}
while passing from first to second line.  The last equality (\ref{QP}) follows from (\ref{Pfinal}). \end{subequations} 
Just to make it clear, Eq. (\ref{QP}) in fact proves $Q=\Phi' P/\mathcal{L}$, regardless of whether a possibility $\Phi'=0$ occurs or not. We have divided by $\Phi'$ only to make the derivation slightly more compact.

By virtue of (\ref{QP}) the equation (\ref{RicPQ}) becomes
\begin{equation}
\Theta_{\pm}'+\Theta_{\pm}^{2}+h=\left(1\pm\frac{\Phi'}{\mathcal{L}}\right)P.
\end{equation}
To finalize the proof we observe that $1\pm\Phi'/\mathcal{L}$ cannot
vanish. Moreover
\begin{equation}\label{finalproof}
u_{\mathtt{top}/\mathtt{bot}}''+h\, u_{\mathtt{top}/\mathtt{bot}}=\left(\Theta_{\mp}'+\Theta_{\mp}^{2}+h\right)u_{\mathtt{top}/\mathtt{bot}}=\left(1\mp\frac{\Phi'}{\mathcal{L}}\right)P\,u_{\mathtt{top}/\mathtt{bot}},
\end{equation}
so that both $u_{\mathtt{top}}$ and $u_{\mathtt{bot}}$ are solutions
of (\ref{Eq.main}) if and only if $P=0$. Consequently, $u = A u_{\mathtt{top}} +B u_{\mathtt{bot}}$, with arbitrary constants $A$ and $B$, is the general solution of (\ref{Eq.main}) if and only if $P=0$. In the final step we make a simple rearrangement
\begin{eqnarray}
P & = & \frac{\Phi\Phi''}{h-\Phi^{2}}+\Phi'\frac{\Phi'\left(h+3\Phi^{2}\right)-\Phi h'}{\left(h-\Phi^{2}\right)^{2}}+\Phi^{2}+h\nonumber \\
 & = & \frac{\Phi}{h-\Phi^{2}}\left[\Phi''-\Phi'\frac{\Phi'\left(h+3\Phi^{2}\right)-\Phi h'}{\Phi\left(\Phi^{2}-h\right)}-\frac{\Phi^{4}-h^{2}}{\Phi}\right].\label{Prearfin}
\end{eqnarray}
to recognize that $P=0$ if and only if (\ref{ELexp}) holds.
\end{proof}

\begin{proof}[Proof of \textup{\fullref{Corollary}{Cor2}}]
By comparison between (\ref{qmp}) and (\ref{Thetapm}) we recognize that $\Theta_-=\Theta_\mathtt{top}$ and $\Theta_+=\Theta_\mathtt{bot}$. Therefore, the first assertion of \fullref{Corollary}{Cor2} follows directly from (\ref{finalproof}) and \fullref{Theorem}{thm-main2}. The second assertion given as Eq. (\ref{PhiSqrt}) follows immediately from multiplication of (\ref{qmp1}) and (\ref{qmp2}).
\end{proof}

\begin{proof}[Proof of \textup{\fullref{Proposition}{prop22}}]
Since, according to the setting framed in \fullref{Fact}{fact}, both $u_1$ and $u_2$ are non-vanishing on $\left]x_-,x_+\right[$, and are solutions of (\ref{Eq.main}), the functions $\Theta_1 = u'_1/u_1$ and  $\Theta_2 = u'_2/u_2$ are solutions of the Ricatti equation (\ref{Ricatti}). Since $\Phi>0$, without loss of generality we can square (\ref{PhiSqrt12}) getting
\begin{equation}
\Phi^{2}=-\Theta_{1}\Theta_{2}.\label{eq:q1q2}
\end{equation}
We first differentiate this equation
\begin{eqnarray}
2\Phi\Phi' & = & -\left(\Theta_{1}'\Theta_{2}+\Theta_{1}\Theta_{2}'\right)=\left(h+\Theta_{1}^{2}\right)\Theta_{2}+\left(h+\Theta_{2}^{2}\right)\Theta_{1}\nonumber \\
 & = & \left(h+\Theta_{1}\Theta_{2}\right)\left(\Theta_{1}+\Theta_{2}\right).
\end{eqnarray}
While replacing $\Theta_{1}'$ and $\Theta_{2}'$ we used the Ricatti equation (\ref{Ricatti}).
With the help of (\ref{eq:q1q2}) we rewrite this result to the form
\begin{equation}
\Phi'=\left(\Theta_{1}+\Theta_{2}\right)\frac{h-\Phi^{2}}{2\Phi}.\label{eq:q1+q2}
\end{equation}
Yet another differentiation, together with the Ricatti equation (\ref{Ricatti}) and (\ref{eq:q1+q2})
leads to
\begin{equation}
\Phi''=\frac{\Phi\Phi'}{h-\Phi^{2}}\frac{d}{dx}\left(\frac{h-\Phi^{2}}{\Phi}\right)-\left(2h+\Theta_{1}^{2}+\Theta_{2}^{2}\right)\frac{h-\Phi^{2}}{2\Phi}.\label{eq:secndder}
\end{equation}
Since
\begin{equation}
\Theta_{1}^{2}+\Theta_{2}^{2}=\left(\Theta_{1}+\Theta_{2}\right)^{2}-2\Theta_{1}\Theta_{2}=\left(\frac{2\Phi\Phi'}{h-\Phi^{2}}\right)^{2}+2\Phi^{2},
\end{equation}
and
\begin{equation}
\frac{\Phi\Phi'}{h-\Phi^{2}}\frac{d}{dx}\frac{h-\Phi^{2}}{\Phi}=\frac{\Phi'\left(h'-2\Phi\Phi'\right)}{h-\Phi^{2}}-\frac{\Phi\left[\Phi'\right]^{2}}{\Phi^{2}},
\end{equation}
Eq. (\ref{eq:secndder}) boils down to (\ref{finalhelpf}) proving the first assertion of \fullref{Proposition}{prop22}.

The second assertion of \fullref{Proposition}{prop22} follows directly from the fact that (\ref{eq:q1q2}) and (\ref{eq:q1+q2}) are equivalent to the Vieta's formulas for roots of a quadratic equation
\begin{equation}\label{quadra}
\Theta^{2}-\frac{2\Phi\Phi'}{h-\Phi^{2}}\Theta-\Phi^{2}=0.
\end{equation}
The roots of (\ref{quadra}) are easily found to be $\Theta_{\mathtt{top}}$ and $\Theta_{\mathtt{bot}}$, given by (\ref{qmp1}) and (\ref{qmp2}) respectively. Finally, Eqs. (\ref{casesS}) follow from integration of $u'_1/u_1=\Theta_{\mathtt{top}/\mathtt{bot}}$ and  $u'_2/u_2=\Theta_{\mathtt{bot}/\mathtt{top}}$ with assumed initial conditions. 
\end{proof}

\begin{proof}[Proof of \textup{\fullref{Lemma}{thm-main}}]
While working with the map (\ref{map}) we can follow two routes. Either treat
the Wronskian as a constant and, if need be, use it to eliminate a
derivative of one solution, or plug its functional form directly into
the map. We will follow the second one. While calculating the Jacobian
matrix of this map we know that the second derivatives $u''_{1}$
and $u''_{2}$ shall appear. Since both $u_{1}$ and $u_{2}$ are
solutions of (\ref{Eq.main}), we will just use this equation to eliminate the second
derivatives. With this single step, and a bit of algebra, we get
\begin{subequations}\label{jacoterms}\begin{equation}
J=\left(\begin{array}{cc}
J_{Xx} & J_{X\varPhi}\\
J_{Yx} & J_{Y\varPhi}
\end{array}\right),
\end{equation}
where:
\begin{equation} 
    J_{Xx}=\frac{\partial X}{\partial x}=\pm \frac{Y^{2}\left(x,\varPhi\right)}{\varPhi^{2}}\left(\varPhi^{2}-h\left(x\right)\right)\left(\left[u'_{1}\left(x\right)\right]^{2}-\varPhi^{2}u_{1}^{2}\left(x\right)\right),
\end{equation}
\begin{equation}
    J_{X\varPhi}=\frac{\partial X}{\partial \varPhi}=\pm 2\frac{Y^{2}\left(x,\varPhi\right)}{\varPhi} u'_{1}\left(x\right)u_{1}\left(x\right),
\end{equation}
\begin{equation}
    J_{Yx}=\frac{\partial Y}{\partial x}= -2\frac{Y^{2}\left(x,\varPhi\right)}{\varPhi}\left(\varPhi^{2}-h\left(x\right)\right)u'_{1}\left(x\right)u_{1}\left(x\right),
\end{equation}
\begin{equation}
    J_{Y\varPhi}=\frac{\partial Y}{\partial \varPhi}= \frac{Y^{2}\left(x,\varPhi\right)}{\varPhi^{2}}\left(\left[u'_{1}\left(x\right)\right]^{2}-\varPhi^{2}u_{1}^{2}\left(x\right)\right).
\end{equation}
Consequently
\begin{equation}
\det J=\pm\frac{\varPhi^{2}-h\left(x\right)}{\left(\varPhi^{2}u_{1}^{2}\left(x\right)+\left[u'_{1}\left(x\right)\right]^{2}\right)^{2}}.
\end{equation}
\end{subequations}
It becomes clear that $\det J\neq0$ in $\mathcal{H}_0$, so  (\ref{map}) is a local diffeomorphism $\mathcal{H}_0\rightarrow \mathcal{H}$ between spaces. We notice the relations between the components of $J$:
\begin{equation}
    J_{Xx}=\pm\left(\varPhi^{2}-h\left(x\right)\right)J_{Y\varPhi},\quad J_{X\varPhi}=\mp \left(\varPhi^{2}-h\left(x\right)\right)^{-1}J_{Y x},
\end{equation}
so that $J_{Xx}J_{X\varPhi}+J_{Yx}J_{Y\varPhi}=0$.

To finalize the proof, we need to check if 
\begin{equation}
g_{h}=J^{T}g_{\mathbb{H}}J=Y^{-2} J^{T}J,
\end{equation}
with $Y$ inside the  diagonal components of $g_{\mathbb{H}}$ replaced by $Y\left(x,\varPhi\right)$
according to (\ref{DiffY}). Alternatively, we can transform differential forms in  $g_{\mathbb{H}}$ as follows
\begin{equation}
dX=J_{Xx}dx+J_{X\varPhi}d\varPhi,\qquad dY=J_{Yx}dx+J_{Y\varPhi}d\varPhi,
\end{equation} 
and replace $Y\mapsto Y\left(x,\varPhi\right)$. Both methods render
\begin{equation}
g_{h}=\frac{1}{Y^{2}\left(x,\varPhi\right)}\left(\begin{array}{cc}
J_{Xx}^{2}+J_{Yx}^{2} & J_{Xx}J_{X\varPhi}+J_{Yx}J_{Y\varPhi}\\
J_{Xx}J_{X\varPhi}+J_{Yx}J_{Y\varPhi} & J_{X\varPhi}^{2}+J_{Y\varPhi}^{2}
\end{array}\right),
\end{equation}
giving the desired result after completing the square
\begin{equation}\label{tocor}
    \left(\left[u'_{1}\left(x\right)\right]^{2}-\varPhi^{2}u_{1}^{2}\left(x\right)\right)^2+\left(2\varPhi u'_{1}\left(x\right)u_{1}\left(x\right)\right)^2 = \frac{\varPhi^{2}}{Y^{2}\left(x,\varPhi\right)},
\end{equation}
followed by a reduction of the fractions. The result is independent of the $\pm$ sign of $X$, as well as of the parameter $X_0$ which neither enters $J$ nor $Y$.
\end{proof}

\begin{proof}[Proof of \textup{\fullref{Corollary}{Cor1}}]
First assertion follows immediately from the sole assumption of \fullref{Lemma}{thm-main}, referring to the fact that both $u_{1}(x)$ and $u_{2}(x)$ are arbitrary linearly independent solutions. The requirement of non-vanishing determinant, defining the group $GL\left(2,\mathbb{R}\right)$, preserves linear independence of all new linear combinations of the solutions.

Nonlinear partial differential equation (\ref{pardif}), in notation from the proof of \fullref{Lemma}{thm-main}, can be rewritten as
\begin{equation}\label{pardif2}
\left(\frac{\varPhi}{h\left(x\right)-\varPhi^{2}}J_{Y x}\right)^{2}+\left(\varPhi J_{Y \varPhi}\right)^{2}=Y^2\left(x,\varPhi\right).
\end{equation}
If we substitute $J_{Y x}$ and $J_{Y \varPhi}$ according to (\ref{jacoterms}), we immediately get the left hand side of (\ref{tocor}) multiplied by $Y^4\left(x,\varPhi\right)/\varPhi^2$. Consequently, Eq. (\ref{pardif2}) becomes equivalent to (\ref{tocor}).
\end{proof}

\section{Details of examples}\label{Sec5} In this supplementary section we collect some details of the presented examples, which are not essential to follow the main body of the paper. In particular, the following claims were provided in Section \ref{subsecEX} without substantiation. In \fullref{Example}{exa1}:
\begin{itemize}
    \item[(1i)] $\Phi_0 (x)$ in Eq. (\ref{Phi0ex1}) as a general solution of (\ref{ELexp}) for $h(x)=0$.
    \item[(1ii)] Solutions (\ref{specsollinear}) of Eq. (\ref{Eq.main})  for $h(x)=0$, stemming from (\ref{Phi0ex1}).
\end{itemize}
In \fullref{Example}{exa2}:
\begin{itemize}
    \item[(2i)] $\Phi_- (x)=\omega$ as a unique constant solution of (\ref{ELexp}) for $h(x)=-\omega^2$.
    \item[(2ii)] Solutions (\ref{solminus}) of Eq. (\ref{Eq.main}) for $h(x)=-\omega^2$, stemming from $\Phi_- (x)=\omega$.
    \item[(2iii)] Injectivity of the map (\ref{regularexp}).
\end{itemize}
In \fullref{Example}{HOexa}:
\begin{itemize}
    \item[(3i)] Range of the solution (\ref{solcot}).
    \item[(3ii)] Injectivity of the map (\ref{XHO}).
    \item[(3iii)] Semicircle image (\ref{circleHOM}).
\end{itemize}
Below we derive the above sometimes very straightforward results one by one.

(1i). In the most trivial scenario $h\left(x\right)=0$, Eq. (\ref{ELexp}) reduces to 
\begin{equation}\label{h0eq}
\Phi''\left(x\right)=3\frac{\left[\Phi'\left(x\right)\right]^{2}}{\Phi\left(x\right)}+\Phi^{3}\left(x\right).
\end{equation}
Upon a substitution $\Phi\left(x\right)=F^{-1/2}\left(x\right)$,
this equation further reduces to a trivial linear equation $F''\left(x\right)=-2$. Its general solution can be represented as $F\left(x\right)=C_{1}^{2}-\left(x+C_{2}\right)^{2}$, leading to the form (\ref{Phi0ex1}). Note that in this way we do not get all solutions of (\ref{h0eq}), since the equation in question is nonlinear. For example,  the same substitution with opposite sign, $\Phi\left(x\right)=-F^{-1/2}\left(x\right)$, renders the same form of $F(x)$.

(1ii).
Utilizing the basic form of the solution (\ref{Phi0ex1}) we find $\Phi_{0}'\left(x\right)=\left(x+C_{2}\right)\Phi_{0}^{3}\left(x\right)$,
$1+\left(x+C_{2}\right)^{2}\Phi_{0}^{2}\left(x\right)=C_{1}^{2}\Phi_{0}^{2}\left(x\right)$
and
\begin{equation}
\sqrt{\Phi_{0}^{4}\left(x\right)+\left[\Phi_{0}'\left(x\right)\right]^{2}}=\Phi_{0}^{2}\left(x\right)\sqrt{1+\left(x+C_{2}\right)^{2}\Phi_{0}^{2}\left(x\right)}=\left|C_{1}\right|\Phi_{0}^{3}\left(x\right).
\end{equation}
Consequently
\begin{eqnarray}
\Theta_{\mathtt{top}/\mathtt{bot},0}\left(x\right) & = & -\frac{\Phi_{0}'\left(x\right)\mp\sqrt{\Phi_{0}^{4}\left(x\right)+\left[\Phi_{0}'\left(x\right)\right]^{2}}}{\Phi_{0}\left(x\right)}\nonumber\\
 & = & \left[\pm\left|C_{1}\right|-\left(x+C_{2}\right)\right]\Phi_{0}^{2}\left(x\right)=\left[\pm\left|C_{1}\right|+\left(x+C_{2}\right)\right]^{-1},
\end{eqnarray}
where we also used $\Phi_{0}^{-2}\left(x\right)=\left|C_{1}\right|^2-\left(x+C_{2}\right)^2$.
The integration in (\ref{specsol}) gives the logarithm, what immediately leads to  (\ref{specsollinear}).

(2i). We first notice that for $\Phi(x)=\mathrm{const}$, Eq. (\ref{ELexp}) becomes
\begin{equation}
   \Phi^4(x)-h^2(x)=0.
\end{equation}
For $h(x)=-\omega^2$ we then have a single valid (real and positive) constant solution equal to $\Phi_-(x)=\omega$. 

(2ii). Since for  $\Phi_-(x)=\omega$ we have $\Phi'_-(x)=0$ and $h(x)-\Phi^2_-(x)=-2\omega^2$, from (\ref{qmp}) we immediately get $\Theta_{\mathtt{top},-}(x)=\omega$ and $\Theta_{\mathtt{bot},-}(x)=-\omega$. Upon integration of a constant function we obtain (\ref{solminus}). 

(2iii). 
 Even though the two dimensional map (\ref{regularexp}) is intuitively injective since with little algebra it was possible to write down its well-behaved inverse (\ref{inverseexp}), standard methods \cite{Injectivity1,Injectivity2} do not confirm that hypothesis because the Jacobian matrix does not meet sufficient criteria (e.g. it is not a P-matrix). However, a technique established in \cite{GaryH1990} for $\mathbb{R}^{2}\rightarrow\mathbb{R}^{2}$ maps would work after being adjusted to the domain in question. While this technique shall also work on convex subsets of $\mathbb{R}^{2}$, instead of tailoring the proof given in \cite{GaryH1990}, we prefer to study injectivity from the definition. Therefore, to establish injectivity of (\ref{regularexp}) we simply ask if equations $X\left(x_a,\varPhi_a\right)=X\left(x_b,\varPhi_b\right)$ and $Y\left(x_a,\varPhi_a\right)=Y\left(x_b,\varPhi_b\right)$ admit nontrivial solutions, i.e. two distinct points $\left(x_a,\varPhi_a\right)$ and $\left(x_b,\varPhi_b\right)$. Finding $e^{-2\omega\left(x_{a}-x_{b}\right)}$ independently from
each equation we get
\begin{equation}
e^{-2\omega\left(x_{a}-x_{b}\right)}=\frac{\varPhi_{b}^{2}-\omega^{2}}{\varPhi_{b}^{2}+\omega^{2}}\frac{\varPhi_{a}^{2}+\omega^{2}}{\varPhi_{a}^{2}-\omega^{2}}=\frac{\varPhi_{b}}{\varPhi_{a}}\frac{\varPhi_{a}^{2}+\omega^{2}}{\varPhi_{b}^{2}+\omega^{2}}.\label{eq:inj}
\end{equation}
The rightmost equation can further be simplified to the form
\begin{equation}
\frac{\varPhi_{b}^{2}-\omega^{2}}{\varPhi_{a}^{2}-\omega^{2}}=\frac{\varPhi_{b}}{\varPhi_{a}}.
\end{equation}
The solutions for $\Phi_{b}$ are
\begin{equation}
\varPhi_{b}=\varPhi_{a},\qquad\textrm{and}\qquad\varPhi_{b}=-\frac{\omega^{2}}{\varPhi_{a}}.\label{eq:sols}
\end{equation}
The second solution is not in $\mathcal{H}$, while for the first
solution Eq. (\ref{eq:inj}) reduces to $e^{-2\omega\left(x_{a}-x_{b}\right)}=1$,
i.e. $x_{a}=x_{b}$. So the map is injective in $\mathcal{H}$ (it
is easy to see it as well is bijective). 

Interestingly, if the map is extended to $\mathbb{R}^{2}\rightarrow\mathbb{R}^{2}$,
then the second solution in (\ref{eq:sols}) is still not valid,
because then
\begin{equation}
e^{-2\omega\left(x_{a}-x_{b}\right)}=-\frac{\omega^{2}}{\varPhi_{a}^{2}}\frac{\varPhi_{a}^{2}+\omega^{2}}{\varPhi_{b}^{2}+\omega^{2}}<0.
\end{equation}
There are no solutions for $x_a-x_b$ and the map in question is injective on $\mathbb{R}^{2}$. 

(3i). The solution (\ref{solcot}) can be rewritten as
\begin{equation}\label{rewrit}
\Phi_{+}\left(x\right)=\omega\sqrt{\frac{\tan\left[\omega\left(x-\bar{x}\right)\right]}{\tan\left[\omega\left(x-\bar{x}-\tilde{x}\right)\right]}}.
\end{equation}
Since the tangent function is, depending on the integer $k$, either increasing (even $k$) or decreasing (odd $k$) in the considered domain (\ref{domainwithk}),
and the phase shift $\tilde{x}$ is positive, the numerator in (\ref{rewrit}) is
always either bigger (even $k$) or smaller (odd $k$) than the denominator. Consequently, either $\Phi_{+}\left(x\right)>\omega$ for even $k$, or $0<\Phi_{+}\left(x\right)<\omega$ for odd $k$.

(3ii).
We first observe that $\zeta\left(x\right)=\cos^{2}\left[\omega\left(x-x_{0}\right)\right]$ must be injective because $Y\left(x,\varPhi\right)$ defined in (\ref{XHO}) depends on $x$ only through $\zeta(x)$. To this end, we need to restrict the domain as (we have already set $k=0$)
\begin{equation}
x\in\left]\bar{x},\bar{x}+\pi/\left(2\omega\right)\right[.
\end{equation}
In this domain we can use the variable $\zeta$ instead of $x$. To study injectivity in full, as before observing that the Jacobian matrix does not fulfill requirements sufficient for injectivity \cite{Injectivity1,Injectivity2}, we look for nontrivial solutions of  $X\left(\zeta_a,\varPhi_a\right)=X\left(\zeta_b,\varPhi_b\right)$ and $Y\left(\zeta_a,\varPhi_a\right)=Y\left(\zeta_b,\varPhi_b\right)$. 
The second equation can be uniquely solved for $\zeta_{b}$
\begin{equation}
\zeta_{b}\left(\varPhi_{a},\varPhi_{b},\zeta_{a}\right)=\frac{\varPhi_{a}^{2}\varPhi_{b}\zeta_{a}+\omega^{2}\left[\varPhi_{b}\left(1-\zeta_{a}\right)-\varPhi_{a}\right]}{\varPhi_{a}\left(\varPhi_{b}^{2}-\omega^{2}\right)}.\label{zetab}
\end{equation}
We consider values of $\left(\varPhi_{a},\varPhi_{b},\zeta_{a}\right)$
for which $0<\zeta_{b}<1$ as otherwise we know there is no solution
to the original problem. Since we know from (\ref{XHO}) that $\textrm{sign}X=\textrm{sign}\left(\omega-\varPhi\right)$,
without loss of generality we can study
\begin{equation}
X^{2}\left(\zeta_{a},\varPhi_{a}\right)-X^{2}\left(\zeta_{b}\left(\varPhi_{a},\varPhi_{b},\zeta_{a}\right),\varPhi_{b}\right)=\frac{\left(\varPhi_{a}-\varPhi_{b}\right)\left(\varPhi_{a}\varPhi_{b}-\omega^{2}\right)}{\omega^{2}\varPhi_{b}\left[\omega^{2}\left(1-\zeta_{a}\right)+\varPhi_{a}^{2}\zeta_{a}\right]}=0.
\end{equation}
We immediately find that the only solution of this equation is $\varPhi_{a}=\varPhi_{b}$, because either $\varPhi_{a},\varPhi_{b}<\omega$ or $\varPhi_{a},\varPhi_{b}>\omega$.
At the end, from (\ref{zetab}) we recover $\zeta_{b}=\zeta_{a}$.

(3iii).
For fixed $\varPhi$ the function $Y\left(\zeta\right)$ is  monotonic
\begin{equation}
\frac{\partial Y}{\partial\zeta}=\varPhi\frac{\omega^{2}-\varPhi^{2}}{\left[\omega^{2}-\zeta\left(\omega^{2}-\varPhi^{2}\right)\right]^{2}},
\end{equation}
since either $\omega<\varPhi$ or $\varPhi<\omega$. We can then invert $Y\left(\zeta,\varPhi\right)$
\begin{equation}
\zeta\left(\varPhi,Y\right)=\frac{\varPhi-Y\omega^{2}}{Y\left(\varPhi^{2}-\omega^{2}\right)},\label{zeta}
\end{equation}
and plug it in (\ref{XHO}), getting
\begin{equation}
X\left(\varPhi,Y\right)=\frac{\textrm{sign}\left(\omega-\varPhi\right)\sqrt{\left(\varPhi-Y\omega^{2}\right)\left(Y\varPhi-1\right)}}{\omega\sqrt{\varPhi}}.\label{XX}
\end{equation}
As already pointed out, the variable $X$ is either positive
or negative, depending on the sign of $\omega-\varPhi$.

Squaring (\ref{XX}) gives an equation of a circle with the center
shifted towards positive values of $Y$
\begin{equation}
X^{2}+\left(\frac{\varPhi^{2}+\omega^{2}}{2\varPhi\omega^{2}}-Y\right)^{2}=\left(\frac{\varPhi^{2}-\omega^{2}}{2\varPhi\omega^{2}}\right)^{2}.
\end{equation}
Since from  (\ref{XHO}) we deduce that $X=0$ for both boundary points $\zeta=0$ and $\zeta=1$, the entire open semicircle is covered (left semicircle if $\omega<\varPhi$ and right semicircle in the opposite case). The summary of all the cases boils down to Eq. (\ref{domM}).

\section{Acknowledgements}
\label{sec:acks}
The author is grateful to Jacek Gulgowski for a helpful discussion. 

\bibliography{aomart}
\bibliographystyle{aomalpha}

\end{document}